%% file: transverse_isotopy_theorem_arxiv.tex
\documentclass[12pt]{amsart}
\usepackage[utf8]{inputenc}
\usepackage{xcolor}
\usepackage[margin=1in]{geometry}
\usepackage{amsmath,amssymb,multicol,cancel,xcolor,graphicx,mathrsfs,amsthm,tikz}
\usepackage[shortlabels]{enumitem}
\usepackage{comment}
\usepackage{thm-restate}
\usepackage{hyperref}

\usepackage{chemfig}
\usepackage[utf8]{inputenc} 
\usepackage[backend=biber,style=alphabetic,]{biblatex}
\hypersetup{
	colorlinks=false, 
	linktoc=all,     
	linkcolor=black,  
}
\usepackage{subcaption}

\usepackage{tcolorbox}

\newcommand{\vertiii}[1]{{\left\vert\kern-0.25ex\left\vert\kern-0.25ex\left\vert #1 
		\right\vert\kern-0.25ex\right\vert\kern-0.25ex\right\vert}}
\newcommand{\R}{\mathbb{R}}

\newcommand{\mc}[1]{\mathcal{#1}}

\newcommand{\bd}{\partial}
\newcommand{\RP}{\mathbb{RP}}

\theoremstyle{definition}
\numberwithin{equation}{section}
\newtheorem{thm}{Theorem}[section]

\newtheorem{cor}[thm]{Corollary}
\newtheorem{defn}[thm]{Definition}

\title{A Transverse Braiding Theorem with Rational Open Book Decomposition}
\author{Ivan So}
\address{Department of Mathematics, Michigan State University, East Lansing, MI 48824}
\email{soivan@msu.edu}

\bibliography{transverse.bib}
\begin{document}
\small
\noindent
\maketitle
\begin{abstract}
	By generalizing the argument of Pavelescu \cite{Pav12}, we show that every transverse link $ K $ in a compact contact 3-manifold can be transversely isotoped to a braid with respect to a rational open book decomposition.
\end{abstract}

\section{Introduction}

In \cite{Ale23}, Alexander showed every link in $ S^3 $ can be presented as a braid closure. The question has been asked in some more general context. Skora \cite{Sko92} proved the same is true for every 3-manifold $ M $ (oriented or not). His idea was to utilize the fact that every 3-manifold admits an (integral) open book decomposition $ (B,\pi) $. Given an open book decomposition $ (B,\pi) $ of $ M $ and a link $ K\subset M$, one can always construct an isotopy such that after the isotopy, $ K $ intersects every page of the open book decomposition transversely and positively.\\
Later, Pavelescu \cite{Pav12} generalized the theorem to the context of contact manifold. Explicitly
\begin{thm}[\cite{Pav12}, Theorem 3.2]
	Suppose $ (B,\pi) $ is an open book decomposition for the 3-manifold $ M $ and $ \xi $ is the contact structure supported by $ (B,\pi) $. Let $ K\subset M $ be a transverse link. Then $ K $ can be transversely isotoped to a braid with respect to $ (B,\pi) $.	
\end{thm}

Theorem 1.1 has numerous application in contact topology. For instance, in defining the invariant for transverse link $t(K)\in HFK^-(-Y,K)$ in \cite{BVV13} and proving the slice-Bennequin inequality for the generalized Rasmussen $ s $-invariant in \cite{MMSW23}, Theorem 1.1 was necessary.\\
However, Theorem 1.1 was established without consideration of rational open book decomposition \cite{BE12}.  The theory of rational open book decomposition  was shown to exhibit a lot of interesting properties other than being a generalized class of decompositions of a 3-manifold. For example, it was shown in \cite{BEV12} that certain cablings of rational open book decomposition have Stein filling and some other properties. Also, based on the cabling results from \cite{BEV12}, Hedden-Plamenevskaya \cite{HP13} was able to conclude some rational surgery over a fibered knot over a contact manifold $ (Y,\xi) $ with nonzero Ozsv\'ath-Szab\'o contact invariant has to be tight.  \\
Building on this, a natural question to ask is: if $ (B,\pi) $ is a \textit{rational} open book decomposition of $ M $, can we always transversely isotope a given link $ K\subset M $ to a braid?\\

We provide a positive answer to this question.

\begin{restatable}{thm}{transverse}
	\label{thm:transverse}
	\textit{Suppose $ (B,\pi) $ is a rational open book decomposition for the compact 3-manifold $ M $ such that $ \xi $ is the contact structure supported by $ (B,\pi) $. Let $ L\subset M $ be a transverse link. Then $ L $ can be transversely isotoped to a braid.}
\end{restatable}

As an immediate corollary, we have the topological version in the sense of \cite{Sko92}.

\begin{cor}
	\textit{Suppose $ (B,\pi) $ is a rational open book decomposition for the 3-manifold $ M $ and $ K\subset M $ a link. Then $ K $ can be isotoped to a braid such that $ K $ intersect each page of $ (B,\pi) $ transversely and positively.}
\end{cor}

The proof of our main theorem follows that of Pavelescu. We show that everything generalize over the rational setting.\\
\textbf{Organization}: In Section 2, we review some background regarding (rational) open book decomposition and contact geometry. Section 3 will be the proof of our main theorem.\\
\textbf{Acknowledgment}: This paper is inspired by the author's work on the slice-Bennequin inequality for the $ \RP^3 $ $ s $-invariant. The author would like to thank Matthew Hedden and John Etnyre for some helpful conversations and their comments on the earlier draft of this paper.


\section{Background}
\subsection{Rational open book decomposition and contact geometry} We first recall a fundamental theorem in contact geometry by Gray in \cite{Gra59}.

\begin{thm}[\cite{Gra59}]
	Let $ \{\xi_T\}_{T\in[0,1]} $ be a family of contact structures on a manifold $ M $ that agrees on the the complement of a compact subset $ C\subset\mathring{M} $, then there exists an isotopy $ \Phi_T:M\to M $ such that 
	\begin{enumerate}[(i)]
		\item 
		$ (\Phi_T)_*\xi_1=\xi_T $;
		\item 
		$ \Phi_T $ is the identity outside an open neighborhood of $ C $.
	\end{enumerate}
	Such family of diffeomorphisms can be defined by the flow of a vector field $ X_T\in\xi_T $ which satisfies
	\begin{align}\label{Gray}
		\iota_{X_T}d\alpha_T=H_T\alpha_T-\frac{d\alpha_T}{dT}
	\end{align}
	with $ \{\alpha_T\} $ the family of contact forms for $ \{\xi_T\} $ and $ H_T $ being a certain function.
\end{thm}

\begin{cor}[\cite{Gra59}]
	Let $ \{\xi_T\}_{T\in[0,1]} $ be a family of contact structures on a manifold $ M $ with respect to the family of contact forms $ \{\alpha_T\}_{T\in[0,1]} $ that differ on a compact set $ C\subset\mathring{M} $, then there exists a vector field $ X_T\in\xi_T $ such that 
	\begin{align*}
		\iota_{X_T}d\alpha_T=H_T\alpha_T-\frac{d\alpha_T}{dT}
	\end{align*}
	for some function $ H_T $.
\end{cor}

Open book decompositions have ties to contact topology in light of the Giroux correspondence \cite{Gir02}. We recall the rational version of open book decomposition defined in \cite{BE12}.

\begin{defn}
	A \textit{rational open book decomposition} for a manifold $ M $ is a pair $ (B,\pi) $ consisting of an oriented link $ B \subset M$ and a fibration $ \pi:M\setminus B\to S^1 $ such that no component of $ \pi^{-1}(\theta) $ is meridional for any $ \theta\in S^1 $.
\end{defn}

In \cite{BEV12}, Baker-Etnyre-Van Horn-Morris proved a generalization of Thurston-Winkelnkemper \cite{TW75} for rational open book decomposition.
\begin{thm}[\cite{BEV12}, Theorem 1.7]\label{BEV}
	Let $ (B,\pi) $ be any rational open book decomposition of $ M $. Then there exists a unique contact structure $ \xi_{(B,\pi)} $ that is supported by $ (B,\pi) $.
\end{thm}

\subsection{Braiding with respect to an open book decomposition}

For the subsequent discussion, we recall some notion defined in \cite{Sko92}. Given a link $ L\subset M $, we can add junctions to it so that we can have the following decomposition.

\begin{defn}[\cite{Sko92}]
	A \textit{piecewise transverse link} (PT link in short) $ L $ in a closed 3-manifold $ Y $ with respect to an open book decomposition $ (B,\pi) $ is a link $ L\subset Y\setminus B $ such that $ L=\cup_i s_i $, where $ s_i $'s are closed segments in which 
	\begin{enumerate}[(1)]
		\item 
		$ \mathring{s}_i\cap\mathring{s}_j =\emptyset$ if $ i\neq j $;
		\item 
		each $ s_i $ is transverse to the pages and intersects each page at most once.
	\end{enumerate}
\end{defn}

To discuss whether a given presentation is braided, we need another notion by Skora

\begin{defn}[\cite{Sko92}]
	A segment in a transverse decomposition of a link is \textit{positive} if the orientation inherited from the link agrees with the transverse orientation of $ (B,\pi) $. We call a segment \textit{negative} otherwise.
\end{defn}

So, from this definition of segments, the game of showing every link can be transversely isotoped to a braid boils down to showing each negative segment of $ L $ can be transversely isotoped to a positive segment.\\
We also mention the following statement by Bennequin, which is a first example of Theorem \ref{thm:transverse}.

\begin{thm}[\cite{Ben83}]\label{Bennequin}
	Any transverse link $ L $ in $ (\R^3,\xi_{\text{std}}) $ is transversely isotopic to a link braided about the $ z $-axis.
\end{thm}

This theorem will be helpful in our proof of Theorem \ref{thm:transverse}.


\section{Proof of the main theorem}
\transverse*

\textit{Proof}: The aim is to show, as in \cite{Pav12}, there exists a family of diffeomorphisms which fix the pages setwise and the flow of the family allows one to push the negative segments to the neighborhood of the binding $ B $. In a neighborhood of the binding, it is contactomorphic to $ (\R^3,\xi_{\text{std}}) $ and apply Theorem \ref{Bennequin}.

We first modify the contact form constructed in \cite{BEV12}. From the information of the open book decomposition, we have
$$M=M_\phi\cup_{\psi}\left(\coprod_{|\bd\Sigma_\phi|}S^1\times D^2\right),$$
with 
\begin{itemize}
	\item 
	$M_\phi$ being the mapping torus of the monodromy $ \phi $ of $ \Sigma $,
	\item 
	$\psi$ the isomorphisms to glue the neighborhood of bindings to $ \Sigma_\phi $.
\end{itemize}
We also use $ (\varphi,(r,\theta)) $ to denote the coordinate on $ S^1\times D^2 $ and $ (s,\theta) $ to denote the coordinates of on a component of $ N(\bd\Sigma_\phi)$, a collar neighborhood of $ \bd\Sigma_\phi $.\\

\textit{Step 1: Construction of $ \alpha_T $}.\\
Let $ \lambda\in\Omega^1(\Sigma) $ such that 
\begin{itemize}
	\item 
	$ d\lambda $ is a volume form on $ \Sigma $ and;
	\item 
	$\lambda=sd\theta$ near each boundary component of $ \Sigma_\phi $.
\end{itemize}

Now, consider $ \lambda_{(t,x)}=t\lambda_x+(1-t)(\phi^*\lambda)_x\in\Omega^1(\Sigma\times I) $ for $ (x,t)\in\Sigma\times[0,1] $. Set 
$$\alpha_{K,T}=\frac{1}{T}\lambda_{(t,x)}+Kdt$$
with $T\in[-1,0) $ and  $K\gg 0 $. Thurston-Weinkelnkemper \cite{TW75} showed such form descends to a contact form on $ \Sigma_\phi $. From here on, we will write $ \alpha_{K,T}=\alpha_T $ as a sufficiently large $ K $ will serve our purpose. Meanwhile, for each $ T\in[-1,0) $, $ \xi_T=\ker\alpha_T $ is a contact structure supported by $ (B,\pi) $, hence the uniqueness from Theorem \ref{BEV} implies the whole family of $ \xi_T $ is isotopic to $ \xi $.\\
\textit{Step 2: Extend $ \alpha_{K,T} $ over the binding.}\\
Explicitly, the map $ \psi $ is given by 
\begin{align*}
	\psi: S^1\times N(\bd D^2)&\to N(\bd\Sigma_\phi)\\
	(\varphi,(r,\theta))&\mapsto (-r,p\varphi+q\theta,-q\varphi+p\theta).
\end{align*}
By $ N(\bd D^2) $, we meant a small annulus. Pulling back the contact form $ \alpha_K $ with $ \psi $ gives
$$\psi^*\alpha_T=(-rp-Kq/T)d\varphi+(-rq+pK/T)d\theta.$$
If an extension of $ \psi^*\alpha_T $ over the entire $ S^1\times D^2 $ exists, it has to be of the form $ f_T(r)d\varphi+g_T(r)d\theta $. So, the existence from Theorem \ref{BEV} boils down to the existence of $ f_T,g_T $ such that
\begin{enumerate}
	\item 
	$ f_Tg_T'-f_T'g_T>0 $ (from the definition of contact form);
	\item 
	near $ \bd (S^1\times D^2) $, $ f_T(r)=-rp-qK/T $ and $ g_T(r)=-rq+pK/T $ (agreement with the pullback we computed);
	\item 
	near the core of $S^1\times D^2$, $ f_T(r)=1 $ and $ g_T(r)=r^2 $.
\end{enumerate}
As in \cite{Pav12}, we can show that such $ f_T,g_T $ exists. Hence there exists a contact structure supported by $ (B,\pi) $. We denote the family of contact forms also by $ \alpha_T $ by abuse of notation.\\

\textit{Step 2: Construction of a $ T $-family of diffeomorphisms $ \{\Phi_T\} $ that fix pages.}\\
The concerned family of diffeomorphisms can be visualized as the flow vector field $ X_T $. From the above steps, we constructed a family of isotopic contact structures. Gray's theorem implies there exists a family of diffeomorphism $ \{\Phi_T\} $, with flow vector $ X_T\in\xi_T $ such that (\ref{Gray}) is satisfied.\\
The requirement of \textit{fixing pages} can be interpreted as $ X_T\in T\Sigma $. Let $ v\in T\Sigma\cap\xi_T $. (\ref{Gray}) applied to $ v $ implies $ d\alpha_T(X_T,v)=0 $. As $ d\alpha_T $ is an area form on $ \xi_T $, this implies $ X_T $ and $ v $ are linearly dependent. So, the family of diffeomorphisms from Gray's theorem is indeed the page fixing family of diffeomorphisms that we are looking for.\\

\textit{Step 3: Singularity analysis of $ X_T $.}\\
Now, we will look for the problematic points of the flow $ X_T $. Note that:
\begin{itemize}
	\item 
	There is no negative elliptic singularities on $ \mathring{\Sigma}_\theta $, as to have elliptic singularity, $ T_e\Sigma_\theta=\xi_e $ at the singularity $ e\in M $.
	\item 
	There may exists hyperbolic singularities. Let $ S_\theta\subset\Sigma_\theta $ be the union of stable submanifold of the hyperbolic singularity and $ \mc{S}=\cup S_\theta $ denote the union of stable submanifolds over pages, which is a 2-CW complex.
\end{itemize}

The set $ \mc{S} $ is problematic for the flow. We need a way out to isotope negative arcs lying in $ \mc{S} $.\\

\textit{Step 4: Wrinkling $ L $ to avoid problematic points on $  \mc{S} $.}\\
To show every link can be transversely braided, it is equivalent to showing all the negative segment can be transversely isotoped to positive ones. We have mentioned in the previous step that outside of $ \mc{S} $, points flow towards a neighborhood of the binding. It will be problematic if an interior point of a negative segment lies on $ \mc{S} $. To deal with this scenario, Pavelescu \cite{Pav12} introduced the notion of wrinkling (See Figure 1 for illustration).\\
She showed we can find a small enough wrinkling such that the whole process is a transverse isotopy. Hence, we are left with negative segments that does not intersect with $ \mc{S} $.\\
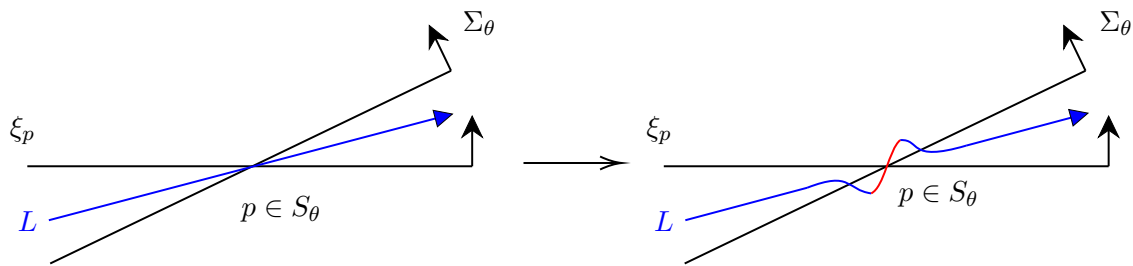
\begin{figure}
	\input{wrinkling}
	\caption{A wrinkling on $ L $ such that the a problematic point is avoided in negative segments.}
\end{figure} 

\textit{Step 5: Applying Bennequin's theorem.} \\
Now, the family of diffeomorphism $ \{\Phi_T\} $ allow us to flow all the negative segments transversely to a neighborhood of the binding. Meanwhile, the neighborhood of the binding is contactomorphic to $ (\R^3,\xi_{\text{std}})/_{z\sim z+1} $ . By Theorem \ref{Bennequin}, negative segments can be transversely isotoped to positive segments. Hence the conclusion of the theorem.\qed

	\printbibliography

\end{document}

%% file: wrinkling.tex
\tikzset{every picture/.style={line width=0.75pt}} 

\begin{tikzpicture}[x=0.75pt,y=0.75pt,yscale=-1,xscale=1]

\draw    (93,121.87) -- (317.36,121.87) ;
\draw    (317.36,121.87) -- (317.36,99.52) ;
\draw [shift={(317.36,96.52)}, rotate = 90] [fill={rgb, 255:red, 0; green, 0; blue, 0 }  ][line width=0.08]  [draw opacity=0] (10.72,-5.15) -- (0,0) -- (10.72,5.15) -- (7.12,0) -- cycle    ;
\draw    (104.54,170.95) -- (306.7,73.65) ;
\draw    (306.7,73.65) -- (297.01,53.52) ;
\draw [shift={(295.71,50.81)}, rotate = 64.3] [fill={rgb, 255:red, 0; green, 0; blue, 0 }  ][line width=0.08]  [draw opacity=0] (10.72,-5.15) -- (0,0) -- (10.72,5.15) -- (7.12,0) -- cycle    ;
\draw [color={rgb, 255:red, 0; green, 0; blue, 255 }  ,draw opacity=1 ]   (103.86,149.11) -- (304.49,96.26) ;
\draw [shift={(307.39,95.5)}, rotate = 165.24] [fill={rgb, 255:red, 0; green, 0; blue, 255 }  ,fill opacity=1 ][line width=0.08]  [draw opacity=0] (8.93,-4.29) -- (0,0) -- (8.93,4.29) -- cycle    ;
\draw    (414,121.87) -- (638.36,121.87) ;
\draw    (638.36,121.87) -- (638.36,99.52) ;
\draw [shift={(638.36,96.52)}, rotate = 90] [fill={rgb, 255:red, 0; green, 0; blue, 0 }  ][line width=0.08]  [draw opacity=0] (10.72,-5.15) -- (0,0) -- (10.72,5.15) -- (7.12,0) -- cycle    ;
\draw    (424.54,170.95) -- (626.7,73.65) ;
\draw    (626.7,73.65) -- (617.01,53.52) ;
\draw [shift={(615.71,50.81)}, rotate = 64.3] [fill={rgb, 255:red, 0; green, 0; blue, 0 }  ][line width=0.08]  [draw opacity=0] (10.72,-5.15) -- (0,0) -- (10.72,5.15) -- (7.12,0) -- cycle    ;
\draw [color={rgb, 255:red, 0; green, 0; blue, 255 }  ,draw opacity=1 ]   (566.5,111.45) -- (625.05,95.84) ;
\draw [shift={(627.95,95.06)}, rotate = 165.07] [fill={rgb, 255:red, 0; green, 0; blue, 255 }  ,fill opacity=1 ][line width=0.08]  [draw opacity=0] (8.93,-4.29) -- (0,0) -- (8.93,4.29) -- cycle    ;
\draw [color={rgb, 255:red, 0; green, 0; blue, 255 }  ,draw opacity=1 ]   (424.86,149.11) -- (486.3,132.72) ;
\draw [color={rgb, 255:red, 0; green, 0; blue, 255 }  ,draw opacity=1 ]   (486.3,132.72) .. controls (508.8,124) and (505.91,133.78) .. (518.58,135.56) ;
\draw [color={rgb, 255:red, 0; green, 0; blue, 255 }  ,draw opacity=1 ]   (533.24,108.67) .. controls (543.91,108) and (536.58,120.22) .. (566.5,111.45) ;
\draw [color={rgb, 255:red, 255; green, 0; blue, 0 }  ,draw opacity=1 ]   (518.58,135.56) .. controls (525.24,132.44) and (527.47,113.33) .. (533.24,108.67) ;
\draw    (343.24,120.22) -- (390.8,120.22) ;
\draw [shift={(392.8,120.22)}, rotate = 180] [color={rgb, 255:red, 0; green, 0; blue, 0 }  ][line width=0.75]    (10.93,-3.29) .. controls (6.95,-1.4) and (3.31,-0.3) .. (0,0) .. controls (3.31,0.3) and (6.95,1.4) .. (10.93,3.29)   ;

\draw (86.7,142.9) node [anchor=north west][inner sep=0.75pt]    {$\textcolor[rgb]{0,0,1}{L}$};
\draw (311.2,41.9) node [anchor=north west][inner sep=0.75pt]    {$\Sigma _{\theta }$};
\draw (82.7,94.9) node [anchor=north west][inner sep=0.75pt]    {$\xi _{p}$};
\draw (407.7,142.9) node [anchor=north west][inner sep=0.75pt]    {$\textcolor[rgb]{0,0,1}{L}$};
\draw (632.2,41.9) node [anchor=north west][inner sep=0.75pt]    {$\Sigma _{\theta }$};
\draw (403.7,94.9) node [anchor=north west][inner sep=0.75pt]    {$\xi _{p}$};
\draw (199.5,135.4) node [anchor=north west][inner sep=0.75pt]    {$p\in S_{\theta }$};
\draw (531,127.9) node [anchor=north west][inner sep=0.75pt]    {$p\in S_{\theta }$};

\end{tikzpicture}

%% file: transverse.bib
@article{Ale23,
  title={A lemma on systems of knotted curves},
  author={Alexander, James W},
  journal={Proceedings of the National Academy of Sciences},
  volume={9},
  number={3},
  pages={93--95},
  year={1923},
  publisher={National Acad Sciences}
}

@article{Ben83,
  title={Entrelacements et {\'e}quations de Pfaff},
  author={Bennequin, Daniel},
  journal={Ast{\'e}risque},
  volume={107},
  pages={87--161},
  year={1983}
}

@article{BE12,
  title={Rational linking and contact geometry},
  author={Baker, Kenneth and Etnyre, John},
  journal={Perspectives in Analysis, Geometry, and Topology: On the Occasion of the 60th Birthday of Oleg Viro},
  pages={19--37},
  year={2012},
  publisher={Springer}
}

@article{BEV12,
  title={Cabling, contact structures and mapping class monoids},
  author={Baker, Kenneth L and Etnyre, John B and Van Horn-Morris, Jeremy},
  journal={Journal of Differential Geometry},
  volume={90},
  number={1},
  pages={1--80},
  year={2012},
  publisher={Lehigh University}
}

@article{BVV13,
  title={On the equivalence of Legendrian and transverse invariants in knot Floer homology},
  author={Baldwin, John A and Vela-Vick, David and V{\'e}rtesi, Vera},
  journal={Geometry \& Topology},
  volume={17},
  number={2},
  pages={925--974},
  year={2013},
  publisher={Mathematical Sciences Publishers}
}

@article{Gra59,
  title={Some global properties of contact structures},
  author={Gray, John W},
  journal={Annals of Mathematics},
  volume={69},
  number={2},
  pages={421--450},
  year={1959},
  publisher={JSTOR}
}

@inproceedings{Gir02,
  title={G\'eom\'etrie de contact: de la dimension trois vers les dimensions sup\'erieures},
  author={Giroux, Emmanuel},
  booktitle={Proceedings of the international congress of mathematicians, ICM 2002},
  pages={405--414},
  year={2002},
  organization={Beijing: Higher Education Press}
}

@article{HP13,
  title={Dehn surgery, rational open books and knot Floer homology},
  author={Hedden, Matthew and Plamenevskaya, Olga},
  journal={Algebraic \& Geometric Topology},
  volume={13},
  number={3},
  pages={1815--1856},
  year={2013},
  publisher={Mathematical Sciences Publishers}
}

@article{MMSW23,
  title={A generalization of Rasmussen’s invariant, with applications to surfaces in some four-manifolds},
  author={Manolescu, Ciprian and Marengon, Marco and Sarkar, Sucharit and Willis, Michael},
  journal={Duke Mathematical Journal},
  volume={172},
  number={2},
  pages={231--311},
  year={2023},
  publisher={Duke University Press}
}

@article{Pav12,
  title={Braiding knots in contact 3-manifolds},
  author={Pavelescu, Elena},
  journal={Pacific journal of mathematics},
  volume={253},
  number={2},
  pages={475--487},
  year={2012},
  publisher={Mathematical Sciences Publishers}
}

@article{Sko92,
  title={Closed braids in 3-manifolds},
  author={Skora, Richard K},
  journal={Mathematische Zeitschrift},
  volume={211},
  pages={173--187},
  year={1992},
  publisher={Springer}
}

@article{TW75,
  title={On the existence of contact forms},
  author={Thurston, William P and Winkelnkemper, Horst E},
  journal={Proceedings of the American Mathematical Society},
  pages={345--347},
  year={1975},
  publisher={JSTOR}
}
